\documentclass[12pt,reqno]{amsart}
\usepackage{fullpage,graphicx}

\newcommand{\jar}{Jarn\'\i k}
\newcommand{\ep}{\varepsilon}
\newcommand{\R}{{\mathbb R}}
\newcommand{\Z}{{\mathbb Z}}

\newtheorem{theorem}{Theorem}
\newtheorem{lemma}[theorem]{Lemma}

\begin{document}

\title{The Limiting Curve of Jarn\'\i k's Polygons}
\author{Greg Martin}
\address{Department of Mathematics\\University of British
Columbia\\Room 121, 1984 Mathematics Road\\Vancouver, BC V6T 1Z2}
\email{gerg@math.ubc.ca}
\maketitle

\section{Introduction}

\begin{figure} 
\includegraphics[width=4.5in,height=3.38in]{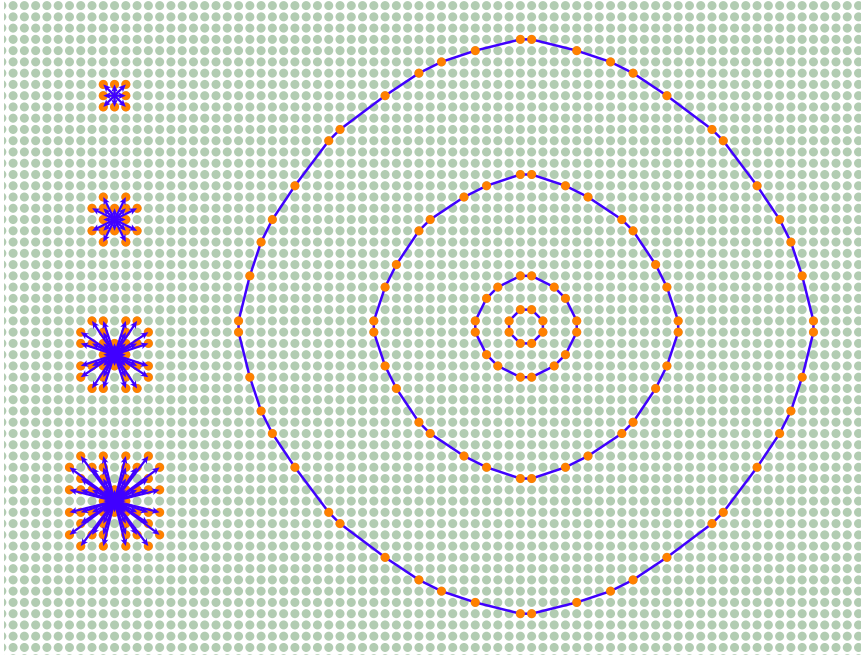}
\caption{\jar\ polygons, right, and their generating sets of vectors, left}
\label{casselman}
\end{figure}

\begin{figure}
\includegraphics[width=2in,height=!]{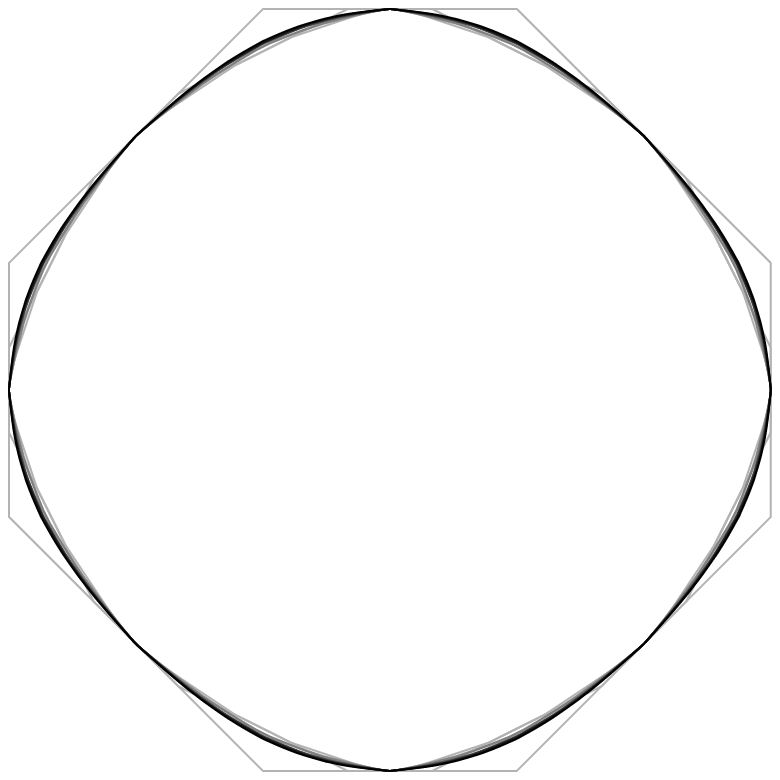}
\hfill
\includegraphics[width=3.46in,height=1.02in]{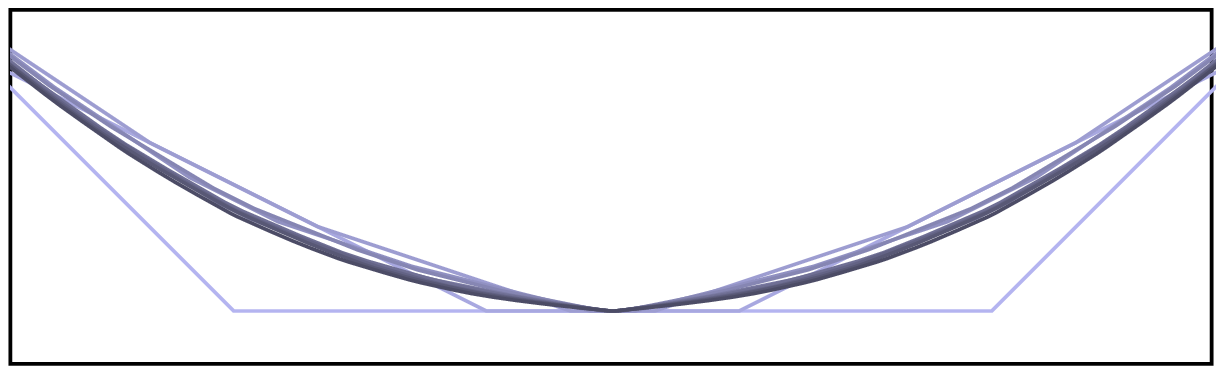}
\caption{Scaled \jar\ polygons superimposed}
\label{superimposed}
\end{figure}

In 1925, \jar\ \cite{jarnik} defined a sequence of convex polygons
for use in constructing curves containing many lattice points relative to
their curvature. Given a positive integer $Q$, let $V_Q$ denote the set of all
primitive integral vectors in the square of side length $2Q$ centered at the
origin, that is,
\begin{equation*}
V_Q = \{ (q,a)\in\Z^2\colon \gcd(q,a)=1,\, \max\{|a|,|q|\}\le Q\}.
\end{equation*}
Then the \jar\ polygon $P_Q$ is the unique (up to translation) convex polygon
whose sides are precisely the vectors in $V_Q$. In other words,
$P_Q$ is the polygon whose vertices can be obtained by starting from an
arbitrary point in $\R^2$ and adding the vectors in $V_Q$ one by one,
traversing those vectors in a counterclockwise direction. For example, the
forty-eight vectors in $V_4$, listed in counterclockwise order, are
\begin{multline}
\qquad\dots,\text{ (1,0), (4,1), (3,1), (2,1), (3,2), (4,3), (1,1),} \\
\text{(3,4), (2,3), (1,2), (1,3), (1,4), (0,1), ($-1$,4), }\dots,\qquad
\label{V4}
\end{multline}
and hence $P_4$ is the tetracontakaioctagon that can be translated to have
vertices at
\begin{multline}
\dots,\text{ ($-1$,0), (0,0), (4,1), (7,2), (9,3), (12,5), (16,8), (17,9),} \\
\text{(20,13), (22,16), (23,18), (24,21), (25,25), (25,26), (24,30), }\dots.
\label{P4}
\end{multline}
These polygons were featured on a recent cover of the {\it Notices of the
American Mathematical Society\/} in connection with an article of Iosevich
\cite{iosevich} and are discussed in further detail in \cite[Chapter
2]{Huxley}. Figure \ref{casselman} shows the four sets of vectors
$V_1$ through $V_4$ and the four polygons $P_1$ through $P_4$ which they
generate.\footnote{Figure \ref{casselman} and the boxed portion of Figure
\ref{superimposed} are a modification of the cover image for the June/July
2001 issue of the {\it Notices of the AMS\/}; they were drawn by Bill
Casselman in Postscript.}

The polygons $P_Q$ have the same eight-fold dihedral symmetry as the unit
square $E=[-1,1]^2$ and, if properly scaled and translated, can be made to pass
through the points $(\pm1,0)$ and $(0,\pm1)$. We denote by $\tilde P_Q$ these
scaled and translated copies of $P_Q$. In a coda to \cite{iosevich}, Casselman
suggests, based on empirical evidence, ``that the scaled polygons [$\{\tilde
P_Q\}$] converge to a somewhat ragged limit curve''. The first several scaled
\jar\ polygons have been superimposed in Figure~\ref{superimposed}, with a
magnified portion shown in the box to the right; the darker polygons correspond to larger values of $Q$. The purpose of this paper is to calculate explicitly the limiting curve of the \jar\ polygons. Indeed, we present several variations on \jar's polygons and calculate the corresponding limit curves, in many cases explicitly and in
other cases parametrically.

If $A$ is a curve in $\R^2$ and $\ep>0$, let $A(\ep)$ denote the
$\ep$-neighborhood of $A$, that is, the set of all points whose distance to
$A$ is less than $\ep$. Given a sequence of curves $A_1, A_2, \dots$, we say
that the curves $\{A_j\}$ converge to $A$ if for every $\ep>0$, there is some
integer $j(\ep)$ such that $A_j$ is contained in $A(\ep)$ for every $j>j(\ep)$.
Our main result, which we prove in Section~2, is the following theorem.

\begin{theorem}
Let $C$ be the curve that contains the graph of the equation
\begin{equation*}
\textstyle y = \frac34x^2 - 1, \quad {-\frac23} \le x \le \frac23
\end{equation*}
and that is invariant under rotation by $\pi/2$ around the origin. Then the
scaled \jar\ polygons $\{\tilde P_Q\}$ converge to $C$.
\end{theorem}

The curve $C$ is infinitely differentiable everywhere except at the four
points $(\pm\frac23,\pm\frac23)$, where it is only twice differentiable.
Although this limiting curve is surprisingly tame, it is the case that the
local ``curvatures'' of the scaled \jar\ polygons oscillate rather than
tending to the corresponding local curvatures of $C$. We describe this
phenomenon in Theorem~\ref{curvature.thm}, the statement and proof of which
appears in Section~6.

\begin{figure}
\includegraphics[width=3in,height=!]{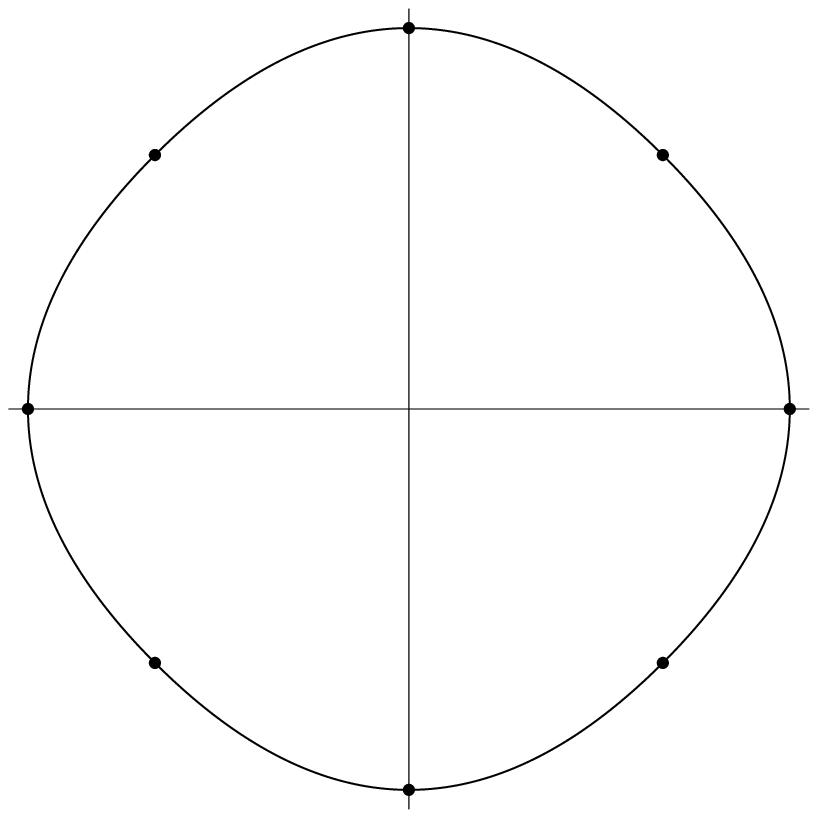}
\hfill
\includegraphics[width=3in,height=!]{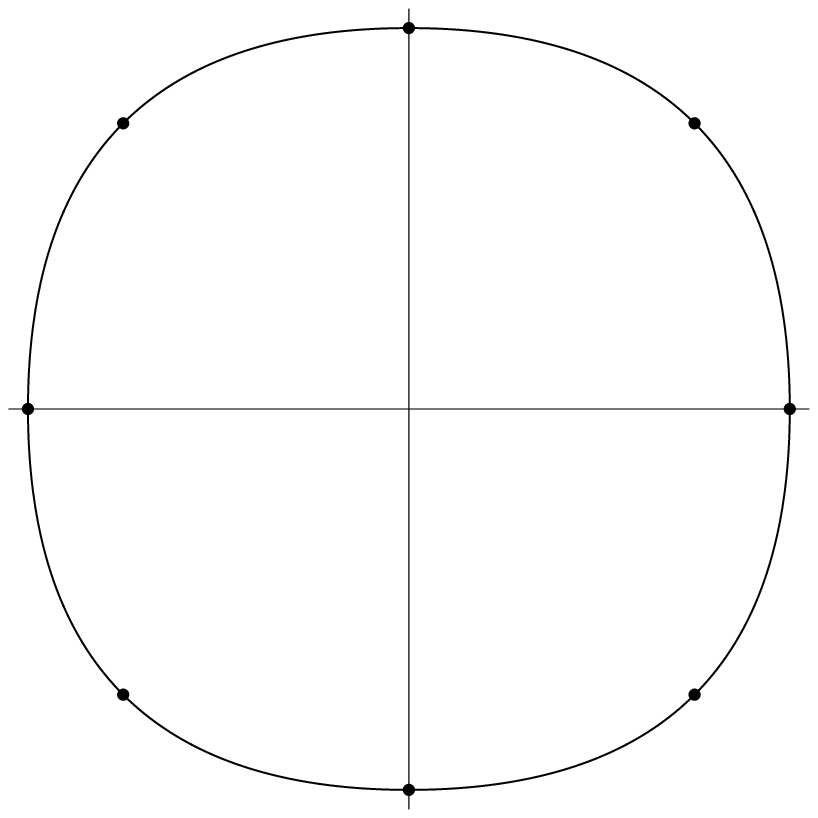}
\begin{picture}(0,0)\small
\put(-330,45){$(\frac23,-\frac23)$}
\put(-287,113){$(1,0)$}
\put(-72,37){$(\frac34,-\frac34)$}
\put(-38,113){$(1,0)$}
\end{picture}
\caption{The curves $C$, left, and $C_1$, right}
\label{CandC1}
\end{figure}

It is interesting to note the relationship between $C$ and the curve $C_1$,
defined as the graph of the equation $\sqrt{1-|x|} + \sqrt{1-|y|} = 1$ (these
two curves are displayed in Figure \ref{CandC1}). Indeed, if $C$ is rotated by
$\frac\pi4$ and then expanded by the factor $\frac3{2\sqrt2}$, so that the
image again passes through the four points $(\pm1,0)$ and $(0,\pm1)$, then the
resulting curve is none other than $C_1$. Vershik \cite{vershik} showed that
$C_1$ is the ``generic'' shape of a convex polygon with lattice point vertices, in the following sense: let $P$ be chosen
at random uniformly from among all convex polygons whose vertices lie in the
set $(\frac1Q\Z)^2 \cap E$. Then, with probability approaching 1 as $Q$ tends
to infinity, the polygon $P$ lies within any prescribed open neighborhood of
$C_1$.

That the curves $C$ and $C_1$ differ only up to rotation and scaling suggests
that rotating the domain demarcating the vectors in $V_Q$ might yield
interesting results. Very generally, given a set $S\subset\R^2$, we may
define the set of vectors
\begin{equation*}
\textstyle V_Q(S) = \{ (a,b)\in\Z^2\colon \gcd(a,b)=1,\, (\frac aQ, \frac bQ)
\in S\}.
\end{equation*}
We also define the corresponding convex polygons $P_Q(S)$ whose sides are
precisely the vectors in $V_Q(S)$, as well as the scaled and translated
versions $\tilde P_Q(S)$ that pass through the four points $(\pm1,0)$ and
$(0,\pm1)$. If we take $S=E$, we recover the vectors $V_Q$ and polygons $P_Q$
in \jar's original definition above. While these definitions make sense for any
set $S$, it seems reasonable in practice to restrict ourselves to sets $S$ that
are star-shaped with respect to the origin and that are the closures of their
interiors; indeed, in this paper we will only consider such sets $S$ centered
at the origin that in addition have the same eight-fold dihedral symmetry as
the unit square. We call the $P_Q(S)$ {\it generalized \jar\ polygons\/}.

Let $D$ denote the ``unit diamond'', namely the square with vertices $(\pm1,0)$
and $(0,\pm1)$. The relationship between $C$ and $C_1$ suggests the following
theorem, which we establish in Section~3.

\begin{theorem}
The scaled generalized \jar\ polygons $\{\tilde P_Q(D)\}$ converge to $C_1$.
\end{theorem}

\noindent Put another way, the generalized \jar\ polygons $P_Q(D)$ are
``generic'' in shape, in the sense of Vershik's theorem.

In this paper we also compute the limit curves corresponding to the $\tilde
P_Q(S)$ for two families of sets $S$, both of which were chosen because they
interpolate between the unit square $E$ and the unit diamond $D$.

\begin{figure}
\includegraphics[width=3in,height=!]{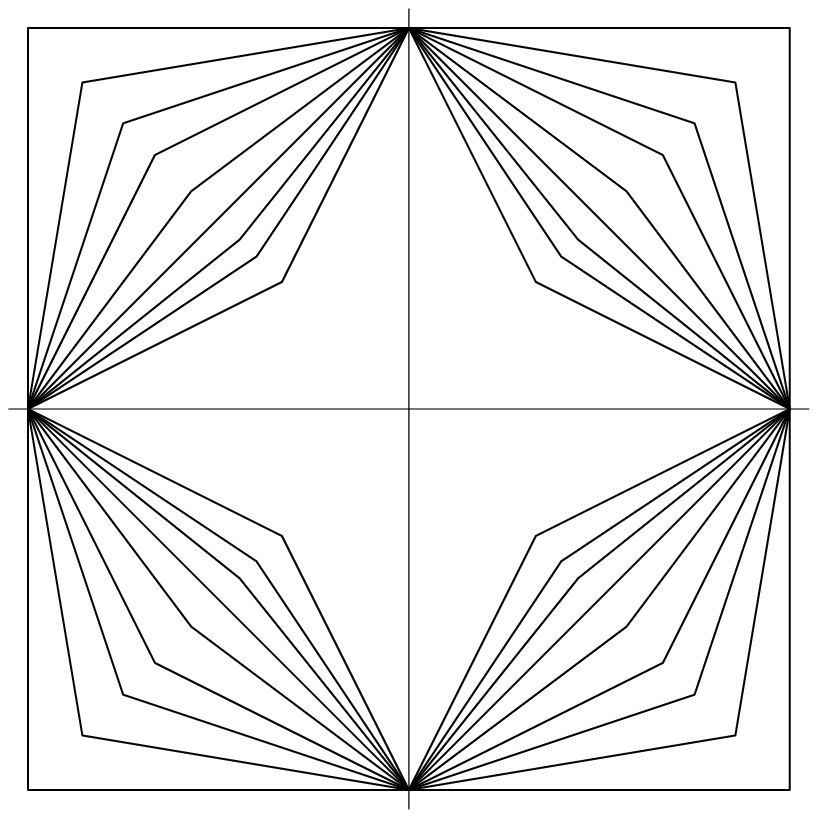}
\hfill
\includegraphics[width=3in,height=!]{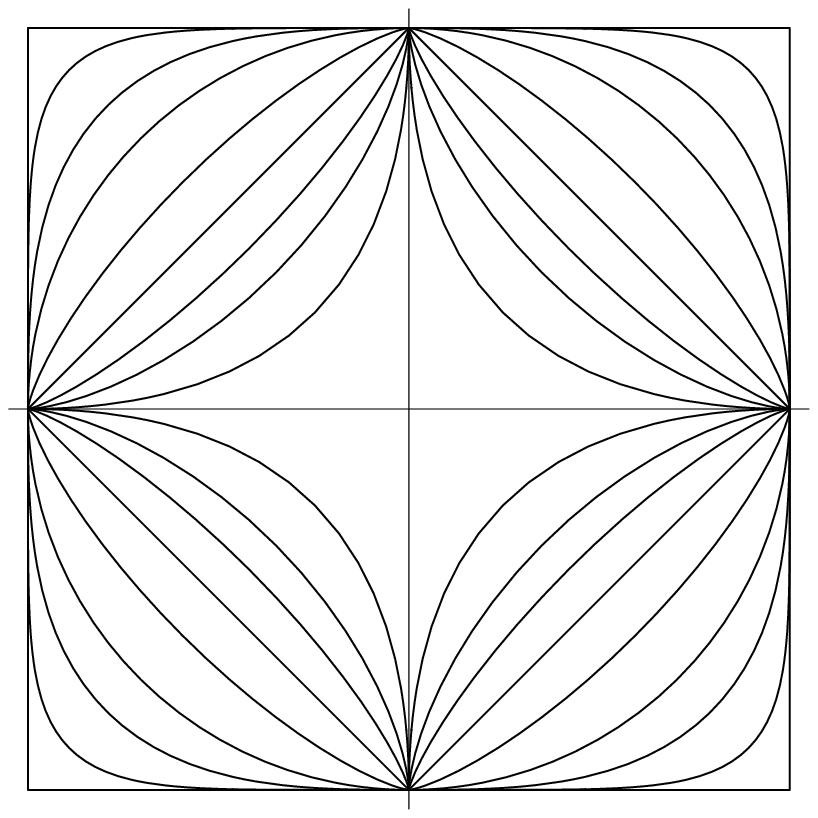}
\caption{The families $O_\delta$, left, and $B_p$, right}
\label{families}
\end{figure}

\begin{itemize}
\item For any positive real number $\delta$, let $O_\delta$ be the octagon
with vertices at $(\pm1,0)$ and $(0,\pm1)$ and the four points
$(\pm\frac{\delta}{1+\delta},\pm\frac{\delta}{1+\delta})$. These octagons
have eight-fold dihedral symmetry, and the slopes of the two edges meeting at
the vertex $(1,0)$ are $\pm\delta$. When $\delta=1$ and as $\delta\to\infty$,
the octagons $O_\delta$ degenerate to the squares $D$ and $E$, respectively.
Figure \ref{families} shows several of these octagons, with $O_{1/2}$ the
innermost and $O_\infty$ the outermost. In Section~4 we calculate the limiting
curves of the polygons $\{\tilde P_Q(O_\delta)\}$ explicitly, and for all
values of $\delta$ these limiting curves are comprised of pieces of parabolas.
\smallskip

\item For any positive real number $p$, let $B_p$ be the set $\{|x|^p + |y|^p
\le 1\}$. When $p\ge1$, the set $B_p$ is simply the closed unit ball in $\R^2$
under the $\ell^p$ metric. Again, when $p=1$ and as $p\to\infty$ we recover the
squares $D$ and $E$. Figure \ref{families} shows several of these sets,
with $B_{1/2}$ the innermost and $B_\infty$ the outermost. (We remark that the
boundary of $B_{1/2}$ is also closely related to Vershik's curve $C_1$.) In
some cases, we can explicitly compute the limiting curves of the $\tilde
P_Q(B_p)$, and in these cases the limiting curves are again piecewise
algebraic. In all cases, we obtain a parametric representation of the limiting
curves and suspect that they are not in general piecewise algebraic. This
family of examples is investigated in Section~5.
\end{itemize}

We note in passing that one can consider the problem of generalizing Vershik's
theorem to domains other than the unit square $E$. Given a set $S\in\R^2$,
choose a polygon $P$ at random uniformly from among all convex polygons whose
vertices lie in the set $(\frac1Q\Z)^2 \cap S$. Is there a curve $V(S)$ such
that, with probability approaching 1 as $Q$ tends to infinity, the polygon $P$
lies within any prescribed open neighborhood of $V(S)$? Vershik's result is
that $V(E)=C_1$; it seems likely that $V(D)$ is the curve $C$ scaled by a
factor of $\frac34$, so that it passes through the points $(\pm\frac12,
\pm\frac12)$. In general, it would be interesting to search for a connection
between these generalized ``Vershik curves'' $V(S)$ and the limiting curves of
familes of generalized \jar\ polygons $\{P_Q(S')\}$ for appropriate sets $S$
and $S'$.

\section{The Original \jar\ Polygons}

Because of the eight-fold symmetry of the \jar\ polygons, we need only
consider the portion of $P_Q$ starting from the edge corresponding to
the vector $(1,0)$ and ending with the edge corresponding to the vector
$(1,1)$; we call this eighth-portion the {\it fundamental arc\/} of $P_Q$, as
the entire polygon $P_Q$ is generated from the fundamental arc under the
action of the dihedral group of order eight. For example, the fundamental arc
of $P_4$ consists of the seven edges defined by the first eight vertices in
equation (\ref{P4}). Given any curve $C$ with eight-fold dihedral symmetry
about the origin (including the scaled \jar\ polygons $\tilde P_Q$ and their
generalizations), we shall also refer to the eighth-portion of the curve lying
in the wedge $\{(x,y)\colon x>0,\, y<-x\}$ as the fundamental arc of~$C$.

For any real number $\lambda\in[0,1]$ we define
$V_Q(\lambda)$ to be the set of all vectors in $V_Q$ with positive coordinates
and slope not exceeding $\lambda$. The sum of all the vectors in
$V_Q(\lambda)$ corresponds to a particular vertex on the fundamental arc of
$P_Q$. If we translate $P_Q$ so that the right-hand endpoint of the edge
corresponding to the vector $(1,0)$ is at the origin, as in (\ref{P4}), then
the coordinates $(X(Q,\lambda),Y(Q,\lambda))$ of this vertex are given by the
formulas
\begin{equation*}
X(Q,\lambda) = \sum_{q\le Q} q \sum_{\substack{a\le\lambda q \\ \gcd(q,a)=1}}
1 \quad\text{and}\quad Y(Q,\lambda) = \sum_{q\le Q}
\sum_{\substack{a\le\lambda q \\ \gcd(q,a)=1}} a.
\end{equation*}
For example, we see from equation (\ref{V4}) that $V_4(\frac1{\sqrt3})$ consists of
the three vectors (4,1), (3,1), and (2,1), and hence
$(X(4,\frac1{\sqrt3}),Y(4,\frac1{\sqrt3}))$ is the vertex (9,3) of $P_4$.

The following asymptotic evaluation of $(X(Q,\lambda)$ and
$Y(Q,\lambda))$ is the key to our calculation.

\begin{lemma}
We have
\begin{equation*}
X(Q,\lambda) = \frac{2\lambda Q^3}{\pi^2} + O(Q^2\log Q) \quad\text{and}\quad
Y(Q,\lambda) = \frac{\lambda^2Q^3}{\pi^2} + O(Q^2\log Q)
\end{equation*}
uniformly for $Q\ge2$ and $0\le\lambda\le1$.
\label{XYlem}
\end{lemma}

{\it Proof:} Recall the definition of the M\"obius mu-function
\begin{equation*}
\mu(n) = \begin{cases}
1, &\text{if } n=1, \\
(-1)^r, &\text{if } n=p_1p_2\dots p_r \text{ where the $p_i$ are distinct
primes,}
\\
0, &\text{if the square of any prime divides $n$.}
\end{cases}
\end{equation*}
The well-known M\"obius inversion formula is based on the characteristic
property of $\mu$
\begin{equation}
\sum_{d\mid n} \mu(d) = \begin{cases} 1, &\text{if }n=1, \\ 0, &\text{if }n>1.
\end{cases}
\label{inversion}
\end{equation}
It is also well known that $\sum_{q\ge1} {\mu(q)}/{q^2} =
{1}/{\zeta(2)} = {6}/{\pi^2}$, where $\zeta$ denotes the Riemann
zeta-function. We shall use the truncated version of this identity
\begin{equation}
\sum_{q\le Q} \frac{\mu(q)}{q^2} = \frac{6}{\pi^2} + O\big( \frac1Q \big),
\label{truncated}
\end{equation}
which follows easily by a trivial estimation of the tail by $\sum_{q>Q}
1/{q^2}$.

Property (\ref{inversion}) allows us to write
\begin{equation*}
X(Q,\lambda) = \sum_{q\le Q} q \sum_{a\le\lambda q} \bigg(
\sum_{d\mid\gcd(q,a)} \mu(d) \bigg) = \sum_{d\le Q} \mu(d) \sum_{\substack{q\le
Q \\ d\mid q}} q \sum_{\substack{a\le\lambda q \\ d\mid a}} 1.
\end{equation*}
Writing $a=bd$ and $q=rd$, we have
\begin{equation}
\begin{split}
X(Q,\lambda) &= \sum_{d\le Q} \mu(d) \sum_{r\le Q/d} dr \sum_{b\le\lambda r} 1
\\
&= \sum_{d\le Q} d\mu(d) \sum_{r\le Q/d} r (\lambda r+O(1)) \\
&= \lambda \sum_{d\le Q} d\mu(d) \sum_{r\le Q/d} r^2 + O\bigg( \sum_{d\le Q} d
\sum_{r\le Q/d} r \bigg) \\
&= \lambda \sum_{d\le Q} d\mu(d) \Big( \frac13 \big(\frac Qd\big)^3 + O\Big(
\big( \frac Qd \big)^2 \Big) \Big) + O\bigg( \sum_{d\le Q} d \big( \frac Qd
\big)^2 \bigg) \\
&= \frac{\lambda Q^3}{3} \sum_{d\le Q} \frac{\mu(d)}{d^2} + O\bigg( Q^2
\sum_{d\le Q} \frac1d \bigg).
\end{split}
\label{middle}
\end{equation}
Equation (\ref{truncated}) now allows us to conclude
\begin{equation*}
X(Q,\lambda) = \frac{\lambda Q^3}{3} \big( \frac{6}{\pi^2} + O\big( \frac1Q
\big) \big) + O(Q^2\log Q) = \frac{2\lambda Q^3}{\pi^2} + O(Q^2\log Q)
\end{equation*}
as claimed.

In the same way we see that
\begin{equation*}
\begin{split}
Y(Q,\lambda) &= \sum_{q\le Q} \sum_{a\le\lambda q} a \bigg(
\sum_{d\mid\gcd(q,a)}  \mu(d) \bigg) \\
&= \sum_{d\le Q} \mu(d) \sum_{\substack{q\le Q \\ d\mid q}}
\sum_{\substack{a\le\lambda q \\ d\mid a}} a \\
&= \sum_{d\le Q} \mu(d) \sum_{r\le Q/d} \sum_{b\le\lambda r} db \\
&= \sum_{d\le Q} d\mu(d) \sum_{r\le Q/d} \big( \frac12(\lambda r)^2 +
O(\lambda r) \big) \\
&= \frac{\lambda^2}{2} \sum_{d\le Q} d\mu(d) \sum_{r\le Q/d} r^2 + O \bigg(
\sum_{d\le Q} d \sum_{r\le Q/d} r \bigg).
\end{split}
\end{equation*}
At this point, a direct comparison to the middle line of equation
(\ref{middle}) yields
\begin{equation*}
Y(Q,\lambda) = \frac{\lambda^2Q^3}{\pi^2} + O(Q^2\log Q)
\end{equation*}
as claimed.\qed\medskip

We can now prove Theorem 1. Define $R(Q) = X(Q,1) + Y(Q,1) - 1/2$, so that
$R(Q) = 3Q^3/\pi^2 + O(Q^2\log Q)$ by Lemma~3. If we translate $P_Q$ so that
the midpoint of the edge corresponding to the vector $(1,0)$ is at the point
$(0,{-R(Q)})$, then the center of $P_Q$ will be at the origin due to the
symmetries of $P_Q$; we then obtain $\tilde P_Q$ by scaling by the factor
$1/R(Q)$. If $(\tilde X(Q,\lambda),\tilde Y(Q,\lambda))$ is the vertex of
$\tilde P_Q$ corresponding to the vertex $(X(Q,\lambda),Y(Q,\lambda))$ of
$P_Q$, then
\begin{equation*}
\begin{split}
\tilde X(Q,\lambda) &= \frac{X(Q,\lambda) +1/2}{R(Q)} =
\frac{2\lambda}{3} + O\big( \frac{\log Q}{Q} \big) \\
\tilde Y(Q,\lambda) &= \frac{Y(Q,\lambda) - R(Q)}{R(Q)} =
\frac{\lambda^2}{3} - 1 + O\big( \frac{\log Q}{Q} \big)
\end{split}
\end{equation*}
using Lemma~\ref{XYlem}. In particular, when $Q$ is large enough, the vertices
on the fundamental arc of $\tilde P_Q$ lie within $\ep/2$ (say) of the arc
parametrized by $(2\lambda/3,\lambda^2/3-1)$ with $0\le\lambda\le1$. This
parametric curve is precisely the arc of the parabola $y = 3x^2/4-1$ from
$x=0$ to $x=2/3$, which is an eighth-portion of the curve $C$. Moreover, the
lengths of the edges of $\tilde P_Q$ are $O(1/Q^2)$, and so every point on the
fundamental arc of $\tilde P_Q$ lies within $\ep$ of $C$ when $Q$ is large
enough. Finally, because of the symmetries of $C$ and the $\tilde P_Q$, we
see that the entire polygon $\tilde P_Q$ lies within an $\ep$-neighborhood of
$C$ when $Q$ is large enough. This establishes Theorem~1.

\section{Polygons Defined by the Unit Diamond}

We begin by looking at the derivation of Lemma \ref{XYlem} from another
viewpoint. For any real number $0\le\lambda\le1$, define $E(\lambda) = \{
(x,y)\in E\colon x>0,\, 0<y\le\lambda x\}$. The inner double sum in the
first line of equation (\ref{middle}) is written as a sum over lattice points
in a large wedge, but we may reinterpret it as a sum over $(\frac1Q\Z)^2 \cap
E(\lambda)$ by writing
\begin{equation*}
\sum_{d\le Q} \mu(d) \sum_{r\le Q/d} dr \sum_{b\le\lambda r} 1 = Q^3
\sum_{d\le Q} \frac{\mu(d)}{d^2} \bigg( \big( \frac dQ \big)^2
\sum_{r\colon 0<dr/Q \le 1} \frac{dr}{Q} \sum_{b\colon 0<db/Q \le \lambda dr/Q}
1 \bigg).
\end{equation*}
Notice that the quantity in parentheses is a Riemann sum approximating the
integral
\begin{equation*}
\int_0^1 x \int_0^{\lambda x} dy\, dx = \iint_{\!\!E(\lambda)} x\, dx\, dy,
\end{equation*}
and in fact (since the integrand $x$ has bounded first derivatives) the error
in making this approximation will be proportional to the mesh size, which is
$O(d/Q)$. Therefore
\begin{align*}
X(Q,\lambda) &= Q^3 \sum_{d\le Q} \frac{\mu(d)}{d^2} \bigg(
\iint_{\!\!E(\lambda)} x\, dx\, dy + O\big( \frac dQ \big) \bigg) \\
&= Q^3 \bigg( \iint_{\!\!E(\lambda)} x\, dx\, dy \bigg) \sum_{d\le Q}
\frac{\mu(d)}{d^2} + O \bigg( Q^2 \sum_{d\le Q} \frac{|\mu(d)|}{d} \bigg) \\
&= \frac{Q^3}{\zeta(2)} \iint_{\!\!E(\lambda)} x\, dx\, dy + O( Q^2\log Q ).
\end{align*}
This is in agreement with Lemma \ref{XYlem}, as $\zeta(2) = \pi^2/6$ and
\begin{equation*}
\iint_{\!\!E(\lambda)} x\, dx\, dy = \int_0^{\lambda} \bigg(
\int_{y/\lambda}^1 x\, dx \bigg) dy = \frac\lambda3.
\end{equation*}
Similar remarks apply to $Y(Q,\lambda)$.

We shall use a similar approach for the generalized \jar\ polygons $P_Q(S)$.
For any real number $\lambda\in[0,1]$ we define $V_Q(S,\lambda)$ to be the set
of all vectors in $V_Q(S)$ with positive coordinates and slope not exceeding
$\lambda$. The sum of all the vectors in $V_Q(S,\lambda)$ corresponds to a
particular vertex on the fundamental arc of $P_Q(S)$. If we translate $P_Q(S)$
so that the right-hand endpoint of the edge corresponding to the vector $(1,0)$
is at the origin, then the coordinates $(X_S(Q,\lambda),Y_S(Q,\lambda))$ of
this vertex are given by the formulas
\begin{equation*}
X_S(Q,\lambda) = \sum_{q\le Q} q \sum_{\substack{(q,a)\in V_Q(S) \\ a\le\lambda
q}} 1 \quad\text{and}\quad Y(Q,\lambda) = \sum_{q\le Q}
\sum_{\substack{(q,a)\in V_Q(S) \\ a\le\lambda q}} a.
\end{equation*}
If we define $S(\lambda) = \{ (x,y)\in S\colon x>0,\, 0<y\le\lambda x\}$, then
the same argument as above allows us to conclude that
\begin{equation}
X_S(Q,\lambda) \sim \frac{Q^3}{\zeta(2)} \iint_{\!\!S(\lambda)} x\, dx\, dy
\quad\text{and}\quad  Y_S(Q,\lambda) \sim \frac{Q^3}{\zeta(2)}
\iint_{\!\!S(\lambda)} y\, dx\, dy.
\label{quickndirty}
\end{equation}
This depends of course on $S$ being a ``reasonable'' set. For the sets $S$ we
shall consider, the asymptotic formulas (\ref{quickndirty}) do in fact hold,
with error terms that are $O(Q^2\log Q)$.

We also use the definitions $R_S(Q) = X_S(Q,1) + Y_S(Q,1) - 1/2$ and
\begin{equation*}
\tilde X_S(Q,\lambda) = \frac{X_S(Q,\lambda)+1/2}{R_S(Q)}, \quad
\tilde Y_S(Q,\lambda) = \frac{Y_S(Q,\lambda)-R_S(Q)}{R_S(Q)},
\end{equation*}
so that $(\tilde X_S(Q,\lambda), \tilde Y_S(Q,\lambda))$ will be the
coordinates of the corresponding vertex on the fundamental arc of the scaled
generalized \jar\ polygon $\tilde P_Q(S)$, the translation and scaling chosen
so that the center of $\tilde P_Q(S)$ is the origin and the points $(\pm1,0)$
and $(0,\pm1)$ are midpoints of edges of $\tilde P_Q(S)$.

We now implement this approach with $S$ equaling the unit diamond $D$ to prove
Theorem 2. Again, due to the symmetries of $C_1$ and the $P_Q(D)$, it suffices
to show that the fundamental arcs of the $P_Q(D)$ tend to the fundamental arc
of $C_1$. Given $0\le\lambda\le1$, the set $D(\lambda)$ is the same as
$\{(x,y)\colon x>0,\, x+y\le1,\, y\le\lambda x\} = \{(x,y)\colon 0<y\le
\lambda/(1+\lambda),\, y/\lambda\le x\le1-y\}$. Therefore
\begin{align*}
\iint_{\!\!D(\lambda)} x\, dx\, dy &= \int_0^{\lambda/(1+\lambda)} \bigg(
\int_{y/\lambda}^{1-y} x\, dx \bigg) dy =
\frac{\lambda(2+\lambda)}{6(1+\lambda)^2}
\\
\iint_{\!\!D(\lambda)} y\, dx\, dy &= \int_0^{\lambda/(1+\lambda)} y \bigg(
\int_{y/\lambda}^{1-y} dx \bigg) dy = \frac{\lambda^2}{6(1+\lambda)^2}.
\end{align*}
Using these evaluations in equation (\ref{quickndirty}), we see that
\begin{equation*}
X_D(Q,\lambda) \sim \frac{Q^3\lambda(2+\lambda)}{\pi^2(1+\lambda)^2}, \quad
Y_D(Q,\lambda) \sim \frac{Q^3\lambda^2}{\pi^2(1+\lambda)^2}
\end{equation*}
as $\zeta(2)=\pi^2/6$. This implies that $R_D(Q) \sim Q^3/\pi^2$, and so
\begin{equation}
\begin{split}
\tilde X_D(Q,\lambda) &= \frac{X_D(Q,\lambda) +1/2}{R_D(Q)} \sim
\frac{\lambda(2+\lambda)}{(1+\lambda)^2} \\
\tilde Y_D(Q,\lambda) &= \frac{Y_D(Q,\lambda) - R_D(Q)}{R_D(Q)} \sim
\frac{-(2\lambda+1)}{(1+\lambda)^2}.
\end{split}
\label{C1.parametric}
\end{equation}
If we set $x=\frac{\lambda(2+\lambda)}{(1+\lambda)^2}$ and
$y=\frac{-(2\lambda+1)}{(1+\lambda)^2}$, it is easy to check that
$\sqrt{1-|x|} + \sqrt{1-|y|} = 1$, and hence the curve parametrized by
$\big( \frac{\lambda(2+\lambda)}{(1+\lambda)^2},
\frac{-(2\lambda+1)}{(1+\lambda)^2} \big)$ with $0\le\lambda\le1$ is precisely
the fundamental arc of $C_1$. This establishes Theorem~2.

\section{Polygons Defined by Octagons}

\begin{figure}
\includegraphics[width=3in,height=!]{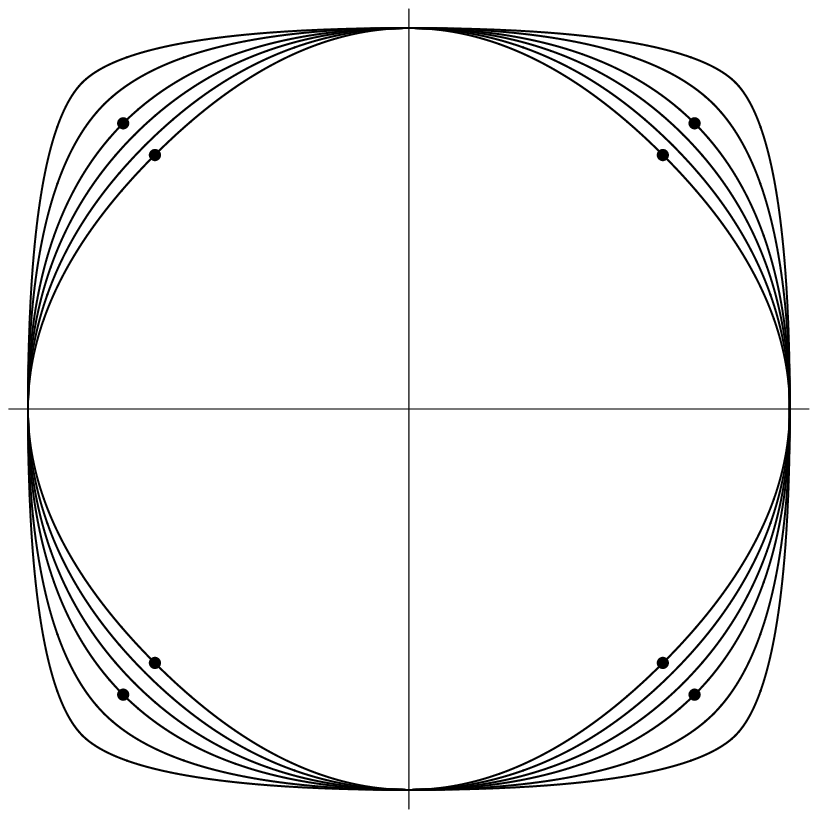}
\hfill
\includegraphics[width=3in,height=!]{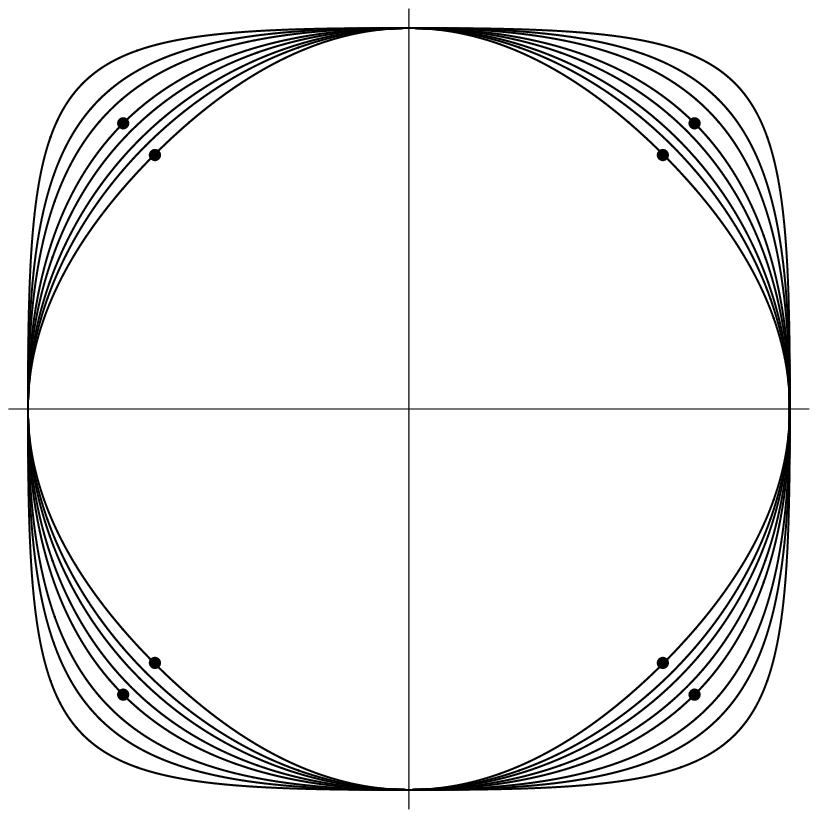}
\caption{The limiting curves $C_\delta$, left, and $C'_p$, right}
\label{limitcurves}
\end{figure}

Recall that for any positive real number $\delta$, we defined $O_\delta$ to be
the octagon with vertices at $(\pm1,0)$ and $(0,\pm1)$ and the four points
$(\pm\frac{\delta}{1+\delta},\pm\frac{\delta}{1+\delta})$. We can use the same
strategy to calculate the limiting curve of the generalized \jar\ polygons
generated from the sets $O_\delta$. Define $C_\delta$ to be the curve with
eight-fold dihedral symmetry whose fundamental arc is
\begin{equation}
\big\{ (x,y)\colon 0\le x\le \tfrac{2\delta+1}{3\delta+1},\,
4\delta(1+\delta)^2(y+1) = (1+3\delta)(\delta x+y+1)^2 \big\}.
\label{delta.parabola}
\end{equation}
This arc is part of a parabola whose axis of symmetry has slope $-\delta$
and whose vertex is $\big( \frac{(1+\delta)^2 (1+2\delta^2)}{(1+3\delta)
(1+\delta^2)^2}, \frac{-(1+2\delta +5\delta^3 +\delta^4
+3\delta^5)}{(1+3\delta) (1+\delta^2)^2} \big)$, as it turns out. The endpoints
of this parabolic arc are $(0,{-1})$ and $(\frac{2\delta+1}{3\delta+1},
-\frac{2\delta+1}{3\delta+1})$. In particular, each $C_\delta$ is piecewise
algebraic, and one can check that each $C_\delta$ is twice differentiable at
the eight symmetry points $(\pm1,0)$, $(0,\pm1)$, and
$(\pm\frac{2\delta+1}{3\delta+1}, \pm\frac{2\delta+1}{3\delta+1})$. Figure
\ref{limitcurves} shows several of these curves, with $C_{1/4}$ the outermost
and $C_\infty$ the innerermost; the points $(\pm\frac23,\pm\frac23)$ and
$(\pm\frac34,\pm\frac34)$, which lie on $C_\infty$ and $C_1$, respectively, are also indicated.

\begin{theorem}
For every positive real number $\delta$, the scaled generalized \jar\ polygons
$\{\tilde P_Q(O_\delta)\}$ converge to $C_\delta$.
\label{octagon.thm}
\end{theorem}

Note that when $\delta=1$, the octagon $O_1$ degenerates to the unit diamond
$D$; in this case, the equation of the parabola in equation
(\ref{delta.parabola}) is equivalent to $4(y+1)=(x+y+1)^2$, which turns out to
be another way to define the fundamental arc of the curve $C_1$. Therefore the
notation $C_\delta$ is consistent with our earlier definition of $C_1$, and
Theorem \ref{octagon.thm} is consistent with Theorem~2. Also, as $\delta$
tends to infinity, the octagons $O_\delta$ converge to the unit square $E$,
while the equation of the parabola in equation (\ref{delta.parabola}) tends to
$4(y+1) = 3x^2$, which is the equation defining the fundamental arc of $C$.
Therefore Theorem \ref{octagon.thm} has Theorem~1 as a limiting case as well.
It can also be checked that $C_\delta$ tends towards the boundary of the unit
square $E$ as $\delta$ decreases to zero.

We prove Theorem \ref{octagon.thm} using the same approach as the proof of
Theorem 2 in the previous section, showing the important steps while
supressing the algebraic details of the computations. Given $0\le\lambda\le1$,
we have by equation (\ref{quickndirty})
\begin{align*}
X_{O_\delta}(Q,\lambda) &\sim \frac{Q^3}{\zeta(2)}
\iint_{\!\!{O_\delta}(\lambda)} x\, dx\, dy \\
&\sim \frac{Q^3}{\zeta(2)} \int_0^{\delta\lambda/(\delta+\lambda)} \bigg(
\int_{y/\lambda}^{1-y/\delta} x\, dx \bigg) dy = \frac{Q^3}{\pi^2}
\frac{\delta\lambda(2\delta+\lambda)}{(\delta+\lambda)^2} \\
Y_{O_\delta}(Q,\lambda) &\sim \frac{Q^3}{\zeta(2)}
\iint_{\!\!{O_\delta}(\lambda)} y\, dx\, dy \\
&\sim \frac{Q^3}{\zeta(2)} \int_0^{\delta\lambda/(\delta+\lambda)} y \bigg(
\int_{y/\lambda}^{1-y/\delta} dx \bigg) dy = \frac{Q^3}{\pi^2}
\frac{\delta^2\lambda^2}{(\delta+\lambda)^2}.
\end{align*}
This implies that $R_{O_\delta}(Q) \sim \frac{Q^3 \delta(3\delta+1)}
{\pi^2(\delta+1)^2}$, and so
\begin{equation*}
\begin{split}
\tilde X_{O_\delta}(Q,\lambda) &= \frac{X_{O_\delta}(Q,\lambda)
+1/2}{R_{O_\delta}(Q)} \sim
\frac{\lambda(2\delta+\lambda)(\delta+1)^2}{(\delta+\lambda)^2(3\delta+1)} \\
\tilde Y_{O_\delta}(Q,\lambda) &= \frac{Y_{O_\delta}(Q,\lambda) -
R_{O_\delta}(Q)}{R_{O_\delta}(Q)} \sim
\frac{\delta\lambda^2(\delta+1)^2}{(\delta+\lambda)^2(3\delta+1)} - 1.
\end{split}
\end{equation*}
If we set $x=\frac{\lambda(2\delta+\lambda)(\delta+1)^2}
{(\delta+\lambda)^2(3\delta+1)}$ and $y=\frac{\delta\lambda^2(\delta+1)^2}
{(\delta+\lambda)^2(3\delta+1)}-1$, one can check that $(x,y)$ satisfies the
polynomial relation (\ref{delta.parabola}), and hence the curve parametrized by
$\big( \frac{\lambda(2\delta+\lambda)(\delta+1)^2}
{(\delta+\lambda)^2(3\delta+1)}, \frac{\delta\lambda^2(\delta+1)^2}
{(\delta+\lambda)^2(3\delta+1)} - 1 \big)$ with $0\le\lambda\le1$ is precisely
the fundamental arc of $C_\delta$. This establishes Theorem~\ref{octagon.thm}.

\section{Polygons Defined by Unit Balls}

Recall that for any positive real number $p$, we defined $B_p$ to be the set
$\{|x|^p + |y|^p \le 1\}$, which we refer to as the ``unit $\ell^p$-ball''
(an abuse of notation when $p<1$). We also need the standard notation
$B(a,b)$ for the Euler beta function $B(a,b) = \int_0^1 t^{a-1}(1-t)^{b-1}\,dt
= \Gamma(a)\Gamma(b)/\Gamma(a+b)$ as well as its relatives, the incomplete
beta function $B_z(a,b) = \int_0^z t^{a-1}(1-t)^{b-1}\,dt$ and the regularized
incomplete beta function $I_z(a,b) = B_z(a,b)/B(a,b)$.

Let $p$ be a positive real number which we regard as fixed. When
$0\le\lambda\le1$, set $\mu=\mu(\lambda)=\frac{\lambda^p}{1+\lambda^p}$, and
define $C'_p$ to be the curve with eight-fold dihedral symmetry whose
fundamental arc is given parametrically by
\begin{equation}
\Big( I_\mu\big( \tfrac1p, 1+\tfrac2p \big) -
\frac{p\lambda(1+\lambda^p)^{-3/p}}{2B(\frac1p,\frac2p)},\,
I_\mu\big( \tfrac2p, 1+\tfrac1p \big) -
\frac{p\lambda^2(1+\lambda^p)^{-3/p}}{B(\frac1p,\frac2p)} - 1 \Big),
\quad0\le\lambda\le1.
\label{Cpp.parametric}
\end{equation}
Figure \ref{limitcurves} shows several of these curves, with $C'_{1/3}$ the
outermost and $C'_\infty$ the innerermost; the points $(\pm\frac23,\pm\frac23)$
and $(\pm\frac34,\pm\frac34)$, which lie on $C'_\infty$ and $C'_1$, respectively, are also indicated.

\begin{theorem}
For every positive real number $p$, the scaled generalized \jar\ polygons
$\{\tilde P_Q(B_p)\}$ converge to $C'_p$.
\label{ball.thm}
\end{theorem}

Although the curves $C'_p$ are in general rather inscrutable, it can be shown
that they are all twice differentiable at the eight points of symmetry and
infinitely differentiable everywhere else. In the special case $p=1$, the ball
$B_1$ is exactly the unit diamond $D$. The parametric representation
(\ref{Cpp.parametric}) of the fundamental arc of $C'_1$ reduces to $\big(
\frac{\lambda(2+\lambda)}{(1+\lambda)^2}, \frac{-(2\lambda+1)}{(1+\lambda)^2}
\big)$, which by equation (\ref{C1.parametric}) is the parametric
representation of the fundamental arc of $C_1$. (Again in this section, we
supress the details of many of our calculations.) Therefore Theorem
\ref{ball.thm} is consistent with Theorem~2.

It can also be shown that as $p$ tends to infinity, the parametric
representation (\ref{Cpp.parametric}) of $C'_p$ approaches
$(2\lambda/3,\lambda^2/3-1)$, which is the parametrization of the fundamental
arc of $C$. Therefore Theorem \ref{ball.thm} has Theorem 1 as a limiting case
as well. Again, it can be checked that $C'_p$ tends towards the boundary
of the unit square $E$ as $p$ decreases to zero.

The case $p=2$ is also special, since the domain $B_2$ (the unit
disk) has complete rotational symmetry. Indeed, the parametric representation
(\ref{Cpp.parametric}) of the fundamental arc of $C'_2$ reduces to $\big(
\frac{\lambda}{\sqrt{1+\lambda^2}}, -\frac{1}{\sqrt{1+\lambda^2}} \big)$,
which is the fundamental arc of the unit circle. Thus the limiting curve of
the generalized \jar\ polygons formed from the vectors in $V_Q(B_2)$ is simply
the unit circle, not surprisingly.

When $p$ is the reciprocal of a positive integer, the regularized incomplete
beta functions are simply indefinite integrals of polynomials. Therefore in
these cases, the parametric representation (\ref{Cpp.parametric}) of the
fundamental arc of $C'_p$ can be written as rational functions of $\lambda^p$.
In particular, these particular curves $C'_p$ are piecewise algebraic, and in
principal one can calculate the algebraic equation defining the fundamental
arc. For example, when $p=1/2$ the parametric representation
(\ref{Cpp.parametric}) of the fundamental arc of $C'_{1/2}$ reduces to
\begin{equation*}
\Big( \frac{\lambda(10+10\lambda^{1/2}+5\lambda+\lambda^{3/2})}
{(1+\lambda^{1/2})^5},\, \frac{-(1+5\lambda^{1/2}+10\lambda+10\lambda^{3/2})}
{(1+\lambda^{1/2})^5} \Big), \quad0\le\lambda\le1,
\end{equation*}
and the coordinates $(x,y)$ of this parametrization satisfy the irreducible
polynomial relation
\begin{multline*}
{-45253} + 86140x - 37030x^2 - 3220x^3 - 765x^4 +
128x^5 - 86140y + 169060xy \\ {}- 80340x^2y -
1940x^3y - 640x^4y - 37030y^2 + 80340xy^2 -
44590x^2y^2 + 1280x^3y^2 \\ {}+ 3220y^3 -
1940xy^3 - 1280x^2y^3 - 765y^4 + 640xy^4 -
128y^5 = 0.
\end{multline*}

We establish Theorem \ref{ball.thm} using our now familiar technique. For any
$0\le\lambda\le1$, equation (\ref{quickndirty}) gives

\begin{align*}
X_{B_p}(Q,\lambda) &\sim \frac{Q^3}{\zeta(2)}
\iint_{\!\!{B_p}(\lambda)} x\, dx\, dy \\
&= \frac{Q^3}{\zeta(2)} \int_0^{\mu^{1/p}} \bigg(
\int_{y/\lambda}^{(1-y^p)^{1/p}} x\, dx \bigg) dy
= \frac{Q^3}{\zeta(2)} \Big( \tfrac1{2p} B_\mu\big(
\tfrac1p, 1+\tfrac2p \big) - \tfrac16 \lambda(1+\lambda^p)^{-3/p} \Big) \\
Y_{B_p}(Q,\lambda) &\sim \frac{Q^3}{\zeta(2)}
\iint_{\!\!{B_p}(\lambda)} y\, dx\, dy \\
&= \frac{Q^3}{\zeta(2)} \int_0^{\mu^{1/p}} y \bigg(
\int_{y/\lambda}^{(1-y^p)^{1/p}} dx \bigg) dy
= \frac{Q^3}{\zeta(2)} \Big( \tfrac1{p} B_\mu\big(
\tfrac2p, 1+\tfrac1p \big) - \tfrac13 \lambda^2(1+\lambda^p)^{-3/p} \Big).
\end{align*}
This implies that
\begin{equation*}
R_{B_p}(Q) \sim \tfrac{Q^3}{\zeta(2)} \tfrac1{3p} B\big( \tfrac1p, \tfrac2p
\big) = \tfrac{Q^3}{\zeta(2)} \tfrac1{2p} B\big( \tfrac1p, 1+\tfrac2p
\big) = \tfrac{Q^3}{\zeta(2)} \tfrac1{p} B\big( \tfrac2p, 1+\tfrac1p
\big),
\end{equation*}
and so after much calculation we see that
\begin{align*}
\tilde X_{B_p}(Q,\lambda) &\sim I_\mu\big( \tfrac1p, 1+\tfrac2p \big) -
\frac{p\lambda}{2(1+\lambda^p)^{3/p}B(\frac1p,\frac2p)} \\
\tilde Y_{B_p}(Q,\lambda) &\sim -\Big( I_{1-\mu}\big( \tfrac1p, 1+\tfrac2p
\big) - \frac{p\lambda^2}{2(1+\lambda^p)^{3/p}B(\frac1p,\frac2p)} \Big),
\end{align*}
which is exactly the parametric definition (\ref{Cpp.parametric}) of the
fundamental arc of $C'_p$. This establishes Theorem~\ref{ball.thm}.

\section{Local Curvatures}

Given a vertex $v$ of any polygon $P$, we define the radius of curvature of
$P$ at $v$ to be the radius of the circle passing through $v$ and its two
neighbors. We quantify the local curvatures of the \jar\ polygons (as
originally defined) in the following way. To each irrational number
$0<\lambda<1$, we associate the unique vertex $v_Q(\lambda)$ on the
fundamental arc of $P_Q$ such that $\lambda$ lies between the slopes of the two
edges adjacent to $v_Q(\lambda)$. We then define $r_Q(\lambda)$ to be the
radius of curvature of $P_Q$ at $v(Q,\lambda)$. This description is not
well-defined for rational numbers $\lambda$, but we can speak of
$r_Q(\lambda^+)$ and $r_Q(\lambda^-)$. For example, from equation (\ref{P4}) we
see that $r_Q(\frac1{\sqrt3})$ is the radius of the circle passing through the
points (7,2), (9,3), and (12,5), which turns out to be $\sqrt{1105/2}$. We
also have $r_Q(\frac12^+)=\sqrt{1105/2}$ but $r_Q(\frac12^-)=5\sqrt{29/2}$.

\begin{figure}
\includegraphics[width=3in,height=!]{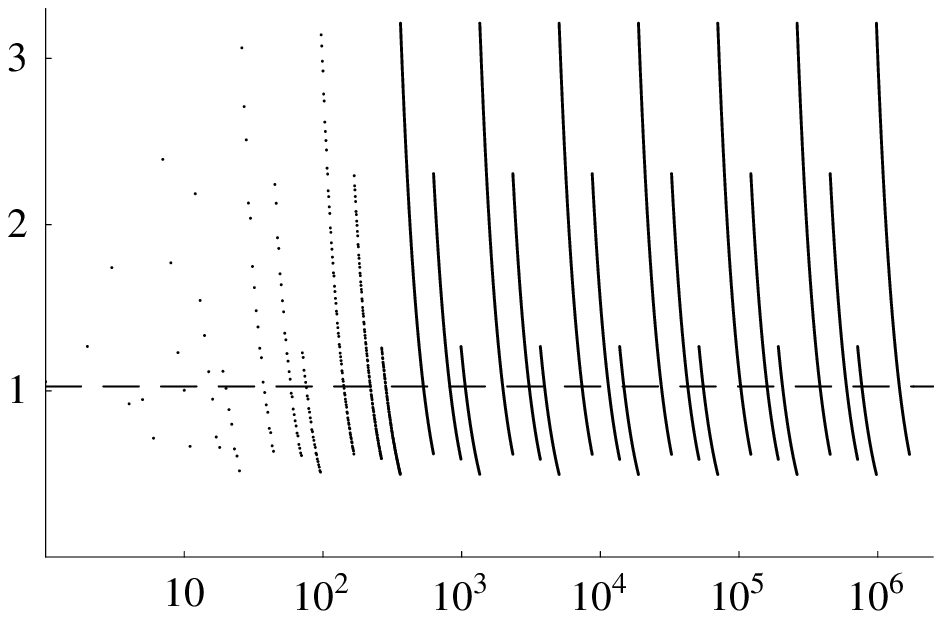}
\hfill
\includegraphics[width=3in,height=!]{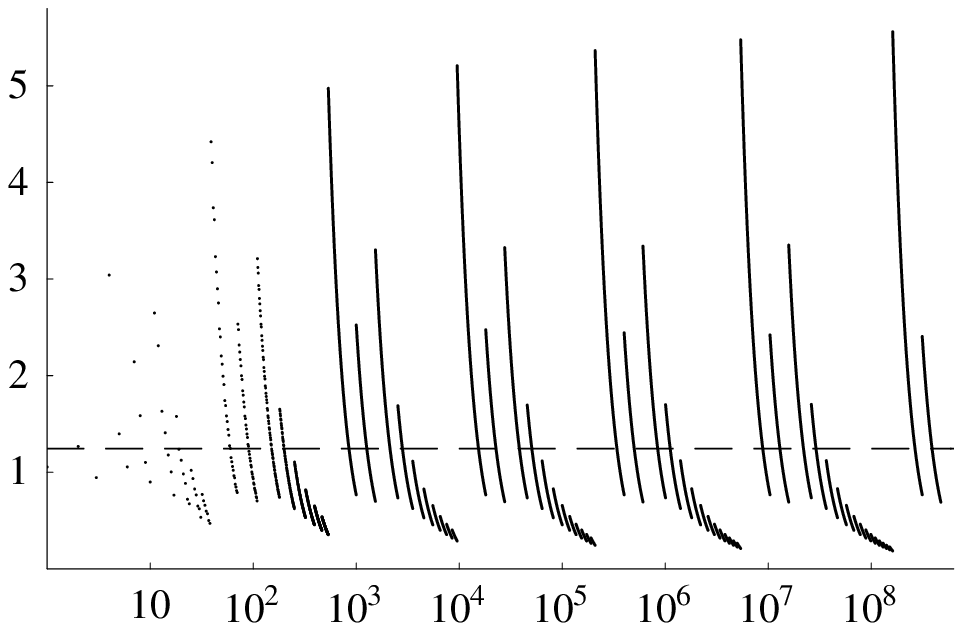}
\caption{The local radii of curvature $\tilde r_Q(\lambda)$ for
$\lambda=\frac1{\sqrt3}$, left, and $\lambda=e-2$, right}
\label{local.curvatures}
\end{figure}

After scaling the \jar\ polygons, the radius of curvature $\tilde
r_Q(\lambda)$ at the corresponding vertex of $\tilde P_Q$ is simply
$r_Q(\lambda)/R(Q)$. It would be tidy if, as $Q$ grew large, the local radii of
curvature $\tilde r_Q(\lambda)$ would converge to the radius of curvature of
the limiting curve $C$ at the corresponding point
$(2\lambda/3,\lambda^2/3-1)$, which turns out to be $\frac23
(1+\lambda^2)^{3/2}$. However, not only does $\{\tilde r_Q(\lambda)\}$ never
converge to $\frac23 (1+\lambda^2)^{3/2}$, but in fact $\{\tilde
r_Q(\lambda)\}$ fails to converge at all for most $\lambda$, and the manner in
which it fails to converge depends upon the diophantine approximation
properties of~$\lambda$. In Figure \ref{local.curvatures} we have plotted these
local radii of curvature $\tilde r_Q(\lambda)$ as functions of $Q$
(represented on the horizontal axis in logarithmic scale) for two interesting
examples of irrational numbers, $\lambda=\frac1{\sqrt3}$ and $\lambda=e-2$,
with the horizontal dashed line indicating the value $\frac23
(1+\lambda^2)^{3/2}$ in each case.

We recall some notation and standard facts from the theory of diophantine
approximation and ontinued fractions. The Farey fractions of order $Q$ are
defined to be the rational numbers in $[0,1]$ with denominator not exceeding
$Q$, listed in increasing order.  For example, the Farey fractions of order 4
are $\{\frac01, \frac14, \frac13, \frac12, \frac23, \frac34, \frac11\}$. There
is a one-to-one correspondence between the Farey fractions $\frac aq$ of order
$Q$ and the vectors $(q,a)$ in the ``fundamental arc'' of $V_Q$, as we see
from equation (\ref{V4}) when $Q=4$. If $\frac{a_1}{q_1}$ and
$\frac{a_2}{q_2}$ are consecutive Farey fractions, it is known (see
\cite[Section 6.1]{NZM}) that $a_2q_1-a_1q_2=1$.

Let $\lambda$ be an irrational number in $(0,1)$ with continued fraction
expansion $[0;b_1,b_2,\dots]$ where the $b_n$ are positive integers. Define
sequences $\{h_n\}$, $\{k_n\}$ of positive integers by
\begin{align*}
h_0=1, \quad h_1=0, \quad h_{n+1} &= b_nh_n+h_{n-1}\;\; (n\ge1) \\
k_0=0, \quad k_1=1, \quad k_{n+1} &= b_nk_n+k_{n-1}\;\; (n\ge1).
\end{align*}
The $\{h_n/k_n\}$ are the {\it convergents\/} to $\lambda$. It is immediate
that for any real number $Q\ge k_2=b_1$, there is a unique index $n\ge2$ and a
unique integer $1\le j\le b_n$ such that $jk_n+k_{n-1} \le Q <
(j+1)k_n+k_{n-1}$. It is known that in the set of Farey fractions of order
$Q$, the number $\lambda$ lies between two fractions whose denominators are
$k_n$ and $jk_n+k_{n-1}$ in this notation. (See \cite[Section 7.5, Problem
5]{NZM}. Fractions of the form $(jh_n+h_{n-1})/(jk_n+k_{n-1})$ with $1\le
j<b_n$ are called the {\it secondary convergents\/} to $\lambda$.)

For example, if $\lambda=\frac1{\sqrt3} = [0;1,1,2,1,2,1,2,\dots]$, then the
sequence of convergents is $\big\{ \frac01, \frac11, \frac12, \frac35, \frac47,
\frac{11}{19}, \frac{15}{26}, \dots \big\}$. Setting $Q=15$, for instance, we
have $k_4=5$, $k_5=7$, $j=1\le2=b_5$, and $1\cdot7+5\le Q<2\cdot7+5$. Hence
the largest (respectively, smallest) rational number with denominator bounded
by 15 that is less than (respectively, greater than) $\frac1{\sqrt3}$ is
$\frac47$ (respectively, $\frac{1\cdot4+3}{1\cdot7+5}=\frac7{12}$), and
therefore the two edges in the fundamental arc of $P_{15}$ between whose
slopes $\frac1{\sqrt3}$ lies correspond to the consecutive vectors (7,4) and
(12,7) of $V_{15}$.

An irrational number $0<\lambda<1$ is {\it badly approximable\/} if
the partial quotients $b_j$ in the continued fraction expansion $\lambda =
[0;b_1,b_2,\dots]$ are bounded, or equivalently if there is a constant
$\delta>0$ such that the inequality $|\lambda-\frac aq| < \frac\delta{q^2}$ has
no solutions. The set of badly approximable irrationals has Lebesgue measure
zero.

\begin{lemma}
For any irrational number $0<\lambda<1$, we have $\liminf (k_n/k_{n+1})
\le(\sqrt5-1)/2$.
\label{Markov.lem}
\end{lemma}

\begin{proof}
Whenever $b_n\ge2$ we have $k_n/k_{n+1} = k_n/(b_nk_n+k_{n-1}) < 1/b_n \le
1/2$. Therefore if infinitely many of the $b_n\ge2$, then $\liminf
(k_n/k_{n+1})\le1/2$. Otherwise, $b_n=1$ for $n$ sufficiently large, so the
$k_n$ eventually satisfy $k_{n+1}=k_n+k_{n-1}$. All solutions to this
recurrence in positive numbers satisfy $k_n \sim c((\sqrt5+1)/2)^n$ for some
constant $c$, hence $\lim (k_n/k_{n+1})=((\sqrt5+1)/2)^{-1} = (\sqrt5-1)/2$.
\end{proof}

We can now describe the limiting behavior of the local radii of
curvature $\tilde r_Q(\lambda)$.

\begin{theorem}
Let $0\le\lambda\le1$.
\begin{enumerate}\renewcommand{\labelenumi}{{\it\alph{enumi}.}}
\item If $\lambda$ is rational, then $\lim_{Q\to\infty} \tilde r_Q(\lambda^+)
= \lim_{Q\to\infty} \tilde r_Q(\lambda^-) = 0$.
\item If $\lambda$ is irrational, then $\limsup_{Q\to\infty} \tilde
r_Q(\lambda)$ lies in the interval
\begin{equation*}
\big[ \tfrac{\pi^2}6 (1+\lambda^2)^{3/2}, \tfrac{\pi^2}3 (1+\lambda^2)^{3/2}
\big].
\end{equation*}
\item $\liminf_{Q\to\infty} \tilde r_Q(\lambda) > 0$ if and only if $\lambda$
is a badly approximable irrational number.
\end{enumerate}
In particular, for almost all $\lambda$, we have $\limsup_{Q\to\infty}
\tilde r_Q(\lambda) > 0$ but $\liminf_{Q\to\infty} \tilde r_Q(\lambda) = 0$.
\label{curvature.thm}
\end{theorem}

\begin{proof}
A straightforward calculation shows that the radius $r$ of the circle passing
through the three points $(x-x_1,y-y_1)$, $(x,y)$, and $(x+x_2,y+y_2)$
satisfies
\begin{equation}
r^2 = \tfrac14(y_1^2+x_1^2)(y_2^2+x_2^2) \big( (y_1+y_2)^2+(x_1+x_2)^2 \big)
(y_2x_1-y_1x_2)^{-2}.
\label{radius}
\end{equation}
To calculate $r_Q(\lambda)$, we take
$(x_1,y_1)=(q_1,a_1)$ and $(x_2,y_2)=(q_2,a_2)$, where
$\lambda$ lies between $\frac{a_1}{q_1}$ and $\frac{a_2}{q_2}$ in the Farey
fractions of order $Q$, and $(x,y)=(X(Q,\lambda),Y(Q,\lambda))$. Since
$\frac{a_1}{q_1}$ and $\frac{a_2}{q_2}$ are consecutive Farey fractions, we
know that $a_2q_1-a_1q_2=1$, and hence the formula (\ref{radius}) simplifies to
\begin{equation}
r_Q(\lambda)^2 = \tfrac14(a_1^2+q_1^2)(a_2^2+q_2^2) \big(
(a_1+a_2)^2+(q_1+q_2)^2 \big).
\label{radius2}
\end{equation}
Now $\frac{a_1}{q_1}$ and $\frac{a_2}{q_2}$ are both approximately $\lambda$,
so substituting $a_1\sim\lambda q_1$ and $a_2\sim\lambda q_2$ into equation
(\ref{radius2}) and simplifying yields
\begin{equation*}
r_Q(\lambda) \sim \tfrac12 q_1 q_2 (q_1+q_2) (1+\lambda^2)^{3/2}.
\end{equation*}
(We record only the main terms for the sake of simplicity. Even a crude
estimate such as $|\frac{a_j}{q_j}-\lambda| = o(\frac1Q)$ would suffice for
our purposes.) Therefore
\begin{equation}
\tilde r_Q(\lambda) = \frac{r_Q(\lambda)}{R(Q)} \sim \frac{q_1 q_2
(q_1+q_2)}{Q^3} \frac{\pi^2 (1+\lambda^2)^{3/2}}6.
\label{rq1q2}
\end{equation}

If $\lambda=\frac aq$ is a rational number, then in using the expression
(\ref{rq1q2}) to calculate $\tilde r_Q(\lambda^+)$, we would have $q_1=q$ for
all $Q\ge q$. In particular, since $q_2\le Q$, the numerator grows only
quadratically with $Q$, and hence $\lim_{Q\to\infty} \tilde r_Q(\lambda^+)=0$.
The same argument shows that $\lim_{Q\to\infty} \tilde r_Q(\lambda^-)=0$, which
proves part (a) of the theorem.

Now suppose that $\lambda=[0;b_1,b_2,\dots]$ is irrational. Then there are
unique positive integers $n$ and $j\le b_n$ such that $jk_n+k_{n-1}\le
Q<(j+1)k_n+k_{n-1}$, where the $k_n$ are the denominators of the convergents to
$\lambda$, in which case $q_1=k_n$ and $q_2=jk_n+k_{n-1}$. Thus the expression
(\ref{rq1q2}) becomes
\begin{equation}
\tilde r_Q(\lambda) \sim \frac{k_n(jk_n+k_{n-1})((j+1)k_n+k_{n-1})}{Q^3}
\frac{\pi^2 (1+\lambda^2)^{3/2}}6.
\label{rks}
\end{equation}
In calculating the lim sup of this expression, we should take $Q$ as small as
possible, that is, $Q=jk_n+k_{n-1}$. If we define $r_n = k_{n-1}/k_n$, then the
expression (\ref{rks}) simplifies to
\begin{equation*}
\tilde r_Q(\lambda) \sim \frac{k_n((j+1)k_n+k_{n-1})}{(jk_n+k_{n-1})^2}
\frac{\pi^2 (1+\lambda^2)^{3/2}}6 = \frac{j+1+r_n}{(j+r_n)^2}
\frac{\pi^2 (1+\lambda^2)^{3/2}}6.
\end{equation*}
The expression $\frac{j+1+r_n}{(j+r_n)^2} = \frac1{j+r_n} + \frac1{(j+r_n)^2}$
is a decreasing function of both $j$ and $r_n$, so in calculating the lim sup
it is best to take $j=1$ and $r_n$ as small as possible. Therefore
\begin{equation*}
\limsup \tilde r_Q(\lambda) = \frac{2+\liminf r_n}{(1+\liminf r_n)^2}
\frac{\pi^2 (1+\lambda^2)^{3/2}}6,
\end{equation*}
and by Lemma \ref{Markov.lem} the first fraction lies in the interval $\big[
\frac{2+(\sqrt5-1)/2}{(1+(\sqrt5-1)/2)^2}, \frac{2+0}{(1+0)^2} \big] = [1,2]$.
This establishes part (b) of the theorem.

Similarly, in calculating the lim sup of the expression (\ref{rks}), we should
take $Q$ as large as possible, that is, $Q=((j+1)k_n+k_{n-1})^-$, in which case
(\ref{rks}) simplifies to
\begin{equation*}
\tilde r_Q(\lambda) \sim \frac{k_n(jk_n+k_{n-1})}{((j+1)k_n+k_{n-1})^2}
\frac{\pi^2 (1+\lambda^2)^{3/2}}6 = \frac{j+r_n}{(j+1+r_n)^2}
\frac{\pi^2 (1+\lambda^2)^{3/2}}6.
\end{equation*}
The expression $\frac{j+r_n}{(j+1+r_n)^2}$ always lies between
$\frac{j+1}{(j+2)^2}$ and $\frac{j}{(j+1)^2}<\frac1j$, which are decreasing
functions of $j$, so in calculating the lim inf it is best to take $j$ as large
as possible, that is, $j=b_n$. If $\lambda$ is badly approximable, so that
$b_n\le B$ for all $n$, then $\liminf \tilde r_Q(\lambda) \ge
\frac{B+1}{(B+2)^2}>0$; if $\lambda$ is not badly approximable, then the $b_n$
are unbounded above and hence
\begin{equation*}
\liminf \tilde r_Q(\lambda) \le \liminf \frac1{b_n} \frac{\pi^2
(1+\lambda^2)^{3/2}}6 = 0.
\end{equation*}
This establishes part (c) of the theorem.
\end{proof}

The two examples in Figure \ref{local.curvatures} illustrate the two
possibilities in part (c) of the theorem. As noted before, the continued
fraction expansion of $\frac1{\sqrt3}$ is $[0;1,1,2,1,2,1,2,\dots]$, and so in
particular $\frac1{\sqrt3}$ is badly approximable since the partial quotients
are bounded above by~2. We can see the repeating groups of a single curve
followed by a pair of curves in the left-hand graph; in particular, the
near-periodicity of the graph implies that the values of $\tilde
r_Q(\frac1{\sqrt3})$ are bounded below. On the other hand, the continued
fraction expansion of $e-2$ is $[0;1,2,1,1,4,1,1,6,1,1,8,\dots]$, and in
particular $e-2$ is not badly approximable since the partial quotients are
unbounded. In the right-hand graph we can see the influence of these partial
quotients (the last full group contains two single curves and a group of 14
curves, corresponding to the string $1,1,14$ in the continued fraction), and
in particular that the lim inf of the values of $\tilde r_Q(e-2)$ is zero.

We remark that the possible values for $\liminf (k_n/k_{n+1})$ in Lemma
\ref{Markov.lem} are closely related to the Markov spectrum. Therefore the
possible values for $\limsup_{Q\to\infty} \tilde r_Q(\lambda)$, as well as the
possible values for $\liminf_{Q\to\infty} \tilde r_Q(\lambda)$ for badly
approximable irrationals $\lambda$, are also related to the Markov spectrum.

\bigskip
{\noindent\small{\it Acknowledgements.} The author acknowledges the support
of the Department of Mathematics of the University of British Columbia and of
the Natural Sciences and Engineering Research Council. The author also
thanks Bill Casselman for pointing out the irregularity of the local curvatures and for contributing the graphics in Figures \ref{casselman} and~\ref{superimposed}.}

\bibliographystyle{amsplain}
\bibliography{jarnik}

\end{document}